\documentstyle[12pt,amssymb,amstex]{article}
\newtheorem{prop}{Proposition}
\newtheorem{thm}{Theorem}
\newtheorem{lm}{Lemma}
\newtheorem{cor}{Corollary}
\newtheorem{conj}{Conjecture}

\title{Correlation for Surfaces of General Type}
\author{Brendan Hassett\thanks{This work was partially supported by a
National Science Foundation Graduate Fellowship}}
\date{\today}
\begin{document}
\maketitle
This paper is still in rough form, and I will appreciate any comments
on the content or the exposition.

\section{Introduction}
The purpose of this paper is to prove the following theorem:
\begin{thm}[Correlation Theorem for Surfaces]
Let $f:X \longrightarrow B$ be a proper morphism of integral
varieties, whose general fiber is an integral surface of general
type.  Then for $n$ sufficiently large, $X^n_B$ admits a dominant
rational map $h$ to a variety $W$ of general type such that the
restriction of $h$ to a general fiber of $f^n$ is generically finite.
\end{thm}
This theorem has a number of geometric and number theoretic consequences
that will be discussed in the final section of this paper.  In particular,
assuming Lang's conjecture on rational points of varieties of general type,
we can prove a uniform bound on the number of rational points on
a surface of general type
not contained in rational or elliptic curves.
\newline \indent This theorem is a special case of the following conjecture
posed by Caporaso, Harris, and Mazur \cite{CHM}:
\begin{conj}[Correlation Conjecture]
Let $f:X \longrightarrow B$ be a proper morphism of integral varieties,
whose general fiber is an integral variety of general type.
\newline Then for $n$ sufficiently large, $X^n_B$ admits a dominant
rational map $h$ to a variety $W$ of general type such that the
restriction of $h$ to a general fiber of $f^n$ is generically finite.
\end{conj}
They prove this conjecture in the case where the general fiber is a curve
of genus $g\geq 2$.  This implies that if Lang's conjectures on
the distribution of rational points on varieties of
general type are true, then there is a uniform
bound on the number of rational points on a curve of genus $g$ defined
over a number field $K$.
The paper \cite{CHM} contains most of the ingredients necessary for
a proof of the general conjecture.  However, at one point the argument
relies heavily on the fact that the fibers of the map are curves:
it invokes the existence of a `nice' class of singular curves, the
stable curves.  For the purposes of this discussion, `nice' means
two things:
\begin{enumerate}
\item Given any proper morphism $f:X\longrightarrow B$
whose generic fiber is a smooth curve of genus $g \geq 2$,
there exists a generically finite
base change $B' \longrightarrow B$ so that the dominating component of
$X \times_B B'$ is birational to a family of stable curves over $B'$.
\item Let $f:X \rightarrow B$ be a family of stable curves, smooth over
the generic point.  Then the fiber products $X_B^n$ are canonical.
\end{enumerate}
For the purpose of generalizing to higher dimensions, we make the following
definitions:
\newline Let $\cal C$ be a class of singular varieties.
\begin{quote}
$\cal C$ is {\em inclusive} if for any proper morphism
$f:X \rightarrow B$ whose generic fiber is a variety of general type,
there is a generically finite base change $B' \rightarrow B$ such that
$X \times_B B'$ is birational to a family $X' \rightarrow B'$ with
fibers in $\cal C$.
\end{quote}
\begin{quote}
$\cal C$ is {\em negligible} if for any family of varieties of general
type ${f:X \rightarrow B}$ with singular fibers belonging to $\cal C$,
the fiber products $X^n_B$ have canonical singularities.
\end{quote}
In a nutshell, the main obstruction to extending the
results of \cite{CHM} is to find a class of higher dimensional singular
varieties which is both negligible and inclusive.
In this paper, we identify a class of surface singularities
which is both inclusive and negligible, and prove that the class
has these properties.  This is the class of
`stable surfaces', surfaces at the
boundary of a compactification of the moduli space of surfaces of
general type.
\newline \indent In the second section of this paper, we describe
these stable surfaces, and try to motivate their definition.  In the
third section, we prove that these stable surfaces actually
form an inclusive class of singularities.
In the fourth section, we prove that stable surfaces
are negligible, i.e. that their fiber products
have canonical singularities.  In the fifth section, we sketch the
proof of the Correlation Theorem for surfaces of general type outlined
in \cite{CHM}.  Finally, we state some consequences of the
Correlation Theorem, assuming various forms of the Lang conjectures.
\vskip.25in
I would like to thank Joe Harris for suggesting this problem, and
Dan Abramovich for his countless comments and corrections.  J\'anos
Koll\'ar also provided useful suggestions for the results
in section four.
\vskip.25in
Thorughout this paper, we work over a field of characteristic zero.
\section{Stable Surfaces}
\indent In this section, we describe the class of stable surfaces
and their singularities.  Most of these results are taken from \cite{K-SB} and
\cite{K}.

\indent
Stable surfaces are defined so that one has a
stable reduction theorem for surfaces, analogous to stable reduction for
curves.  For motivation, we first review the curve case.
Let $f:X \rightarrow \Delta$ be a flat family of curves of genus $g \geq 2$
over a disc.  Assume that the fibers of the
family are smooth, except for the fiber $f^{-1}(0)$ which may be singular.
By Mumford semistable reduction (\cite{KKMS}), there is a finite base change
$$\tilde{\Delta} \longrightarrow \Delta$$
ramified over $0$, and a resolution of singularities of the base-changed family
$$d: Y \longrightarrow \tilde{X}$$
such that the fibers of the composed map
$$F=d \circ \tilde{f}:Y \longrightarrow \tilde{\Delta}$$
are reduced normal crossings divisors.  This semistable reduction is not
unique, as we can always blow up $Y$ to get a `different' semistable reduction.
These semistable reductions are all birational, and we can take a
`canonical model' $\cal Y$ of the surface $Y$ by using the relative
pluricanonical differentials to map $Y$ birationally into projective space.
$\cal Y$ is the image of this map, and  $\cal Y \rightarrow \tilde{\Delta}$ is
called the stable reduction of our original family.  Moreover, the birational
map $Y \rightarrow {\cal Y}$ can be
described quite explicity.  It is the morphism that blows down all the $-1$
and $-2$ curves on $Y$.  On the fibers, this corresponds to blowing
down smooth rational components meeting the rest of the fiber in
one or two points.
{}From this, we see that the fibers of ${\cal Y} \rightarrow \tilde{\Delta}$
are just stable curves.
\newline \indent For higher dimensional varieties, we can try to mimic the same
procedure.  We can still apply semistable reduction to obtain the family
$Y \longrightarrow \tilde{\Delta}$, but this reduction is not unique.
The problem is that it is not generally known how to obtain a canonical model
$\cal Y$ for the birational equivalence class of $Y$.  This canonical model
$\cal Y$ would be our stable reduction, if it were well defined.  In the
case of families of surfaces where $Y$ is a threefold,
we can use the minimal model program to
construct the canonical model of a semistable family of surfaces (cf
\cite{Ka}).
The total space of our stable family will then have canonical singularities,
and the singularities of the fibers of the family can then be described.
\newline \indent Now we introduce the formal definitions.
By definition, a variety $S$ is said to be $\Bbb Q$-Gorenstein if
$\omega_S^{[k]}$ is locally free for some $k$.
$\omega_S^{[k]}$ denotes the reflexive hull
(i.e. the double dual) of the $k$th power of the dualizing
sheaf.  For a $\Bbb Q$-Gorenstein singularity, the smallest such $k$
is called the index of the singularity.
A surface is {\em semi-smooth} if it has only the following singularities:
\begin{enumerate}
\item { $2$-fold normal crossings with equation $x^2=y^2$}
\item { pinch points with equation $x^2=zy^2$ }
\end{enumerate}
A {\em good semi-resolution} resolution of $S$ is a proper map
$g:T \longrightarrow S$ satisfying the following properties
\begin{enumerate}
\item { $T$ is semi-smooth}
\item { $g$ is an isomorphism in the complement of a codimension
two subscheme of $T$}
\item { No component of the double curve of $T$ is exceptional for $g$.}
\item {The components of the double curve of $T$ and the exceptional locus
of $S$ are smooth, and meet transversally.}
\end{enumerate}

A surface $S$ is said to have semi-log-canonical singularities if
\begin{enumerate}
\item $S$ is Cohen-Macaulay and $\Bbb Q$-Gorenstein with index $k$
\item $S$ is semi-smooth in codimension one
\item The discrepancies of a good semi-smooth resolution of $S$ are all
greater than or equal to $-1$
(i.e. $\omega_T^k = g^*\omega^{[k]}_S(ka_1 E_1 +...+ ka_n E_N)$ where $a_i
\geq-1$)
\end{enumerate}
In \cite{K-SB} a complete classification of semi-log-canonical singularities
is given.
\newline \indent The relevance of these definitions comes from a
result proved in the same paper

\begin{thm}
Let $f: X \rightarrow \Delta$ be a family of surfaces over the
disc.  Then the following are equivalent:
\begin{enumerate}
\item The general fiber has rational double points, and the central fiber
has semi-log-canonical singularities.
\item For any base change $\tilde{\Delta} \rightarrow \Delta$, the base-changed
family
$$\tilde f: \tilde X \longrightarrow \tilde{\Delta}$$ has canonical
singularities.
\end{enumerate}
In fact, if $X \rightarrow \Delta$ has a semistable resolution
of singularities, then $X$ is canonical iff the general fiber has rational
double points and the central fiber has semi-log-canonical singularities.
\end{thm}
In particular, this means that the `bad' fibers in a stable
reduction of surfaces have only semi-log-canonical singularities.
For the sake of this discussion, surfaces with only rational double
point singularities are `good' fibers.  This is reasonable, because
we would like the canonical model of a smooth surface of general type
to be `good'.
This motivates the definition:
\begin{quote}
A surface $S$ is {\em stable} if $S$ has semi-log-canonical
singularities, and for some sufficiently large $k$ $\omega_S^{[k]}$ is locally
free and ample.
\end{quote}
Note that a smooth surface of general type is not stable
if it contains $-1$ or $-2$ curves, but its canonical model will be stable.
\begin{quote}
A family of stable surfaces is defined to be a proper flat morphism
$\cal S \rightarrow B$ whose fibers are stable surfaces,
with the property that taking reflexive powers of the relative
dualizing sheaf commutes with restricting to a fiber:
$$\omega_{\cal S/B}^{[k]}|{\cal S}_b = \omega^{[k]}_{{\cal S}_b}$$
\end{quote}
In particular, the reflexive powers of the relative dualizing sheaf are flat.
This additional condition is necessary to guarantee that the  moduli space
in the next section is separated.  Note that we can define
$$K_S^2 = {1\over {k^2}}\#(\omega_S^{[k]},\omega_S^{[k]})$$
for any stable surface $S$, and that this number is constant
in families.  We also have the invariant $\chi_S=\chi(\cal O_S)$, which is
also constant in families.
Finally, stable surfaces are analogous to stable curves in
one more important sense:
\begin{thm}
A stable surface has a finite automorphism group.
\end{thm}
The essence of the proof is easy to grasp.  Let $S$ be stable, and let
$\tilde S$ be its normalization.  Let $\Delta$ be the double curve on
$\tilde S$.
The pair $(\tilde S,\Delta)$, is log-canonical (see \cite{K-SB}).
Therefore, each component of $(\tilde S,\Delta)$ is of log-general type, and
has a finite automorphism group by \cite{I}.

\section{Stable Surface Singularities are Inclusive}
\indent To show that the class of stable surfaces are inclusive,
we need to invoke the existence of a proper coarse moduli space
$\bar{\cal M}_{\chi,K^2}$ for the stable surfaces with invariants $\chi$ and
$K^2$.  We also need a finite covering $\phi:\Omega \rightarrow \bar{\cal
M}_{\chi,K^2}$ of
the moduli space that admits a tautological family $\cal S \rightarrow \Omega$:
 $$\begin{array}{ccc}
\cal T &	&	\\
\downarrow	&	& 	\\
\Omega  &   \stackrel{\phi}{\rightarrow} & \bar{\cal M}_{\chi,K^2}\\
\end{array}$$
We have the following theorem:
\begin{thm}
For smoothable stable surfaces, there exists a coarse moduli space $\bar{\cal
M}_{\chi,K^2}$ with these properties.
\end{thm}
By definition, a stable surface is {\em smoothable} if it is contained in a
family of stable surfaces with $\Bbb Q$-Gorenstein total space, such that the
general member has only rational double points.
The proof of this theorem is scattered throughout the literature.  The proof
that the moduli space exists as a separated algebraic space is contained
in \cite{K-SB} \S 5.  This relies on the properties of semi-log-canonical
singularities and the finite automorphism theorem.
The proof that the moduli space has a functorial
semipositive polarization is contained in \cite{K} \S 5.  This
paper also has a general argument for the existence of a
finite covering of the moduli space possessing a tautological
family (see also \cite{CHM} \S 5.1).
The proof that the moduli space is of finite type for a
given pair of invariants (and thus proper and projective by \cite{K})
is contained in \cite{A}.
\newline \indent Using this moduli space, we can prove that the class of
stable surface singularities is inclusive.
\begin{prop}
Let $f:X \rightarrow B$ be a proper morphism of integral
varieties.  Assume that the general fiber of this map is a
smooth surface of general type.  Then there exists a generically
finite base change $B' \rightarrow B$ such that the dominating
component of the
fiber product $X \times_B B'$ is birational to a family of stable
surfaces over $B'$.
\end{prop}
The proof of this follows the proof
of the analogous result in \cite{CHM} \S 5.2 t quite closely.
Since $\bar{\cal M}$ is a coarse moduli space,
there is an induced rational map $B \rightarrow \bar{\cal M}$.
Let $B_1$ be the closed graph of this map, $\Sigma_1$
its image in $\bar{\cal M}$, and $X_1\rightarrow B_1$
the dominating component
of $X\times_B B_1$.  Since there is no tautological family
on $\bar{\cal M}$, we do not have a family of stable surfaces defined
over $\Sigma_1$.  However, we do have such a family over
$\Sigma_2=\phi^{-1}(\Sigma_1)\subset \Omega$, which we denote
$\cal T_2$.  So we let $B_2=B_1 \times_{\Sigma_1} \Sigma_2$,
$\mu:B_2 \rightarrow \Sigma_2$ the projection, and $X_2 \rightarrow B_2$ the
main component of $X_1 \times_{B_1} B_2$.  The family $\cal T_2$
of stable surfaces pulls back to a family
$Y_2=\cal T_2 \times_{\Sigma_2}B_2 \rightarrow B_2$.  For general $b\in B_2$,
$(X_2)_b$ is birational to $(Y_2)_b$ by construction.

\indent This is almost enough to prove
that $X_2$ is birational to the family of stable surfaces $Y_2$.  We
just need one further finite base change $B_3 \rightarrow B_2$
to `straighten out' $X_2$.  We have to get rid of isotrivial subfamilies that
cannot be represented by pull backs of the tautological
family on the moduli space.  We describe the fiber of this base change over
the generic point $b\in B_2$.  Set theoretically, it will
correspond to equivalence classes of birational morphisms

$$ \{ \psi | \psi: (X_2)_b \rightarrow(\cal T_2)_{\mu(b)}\}$$

Two maps are equivalent if they induce the identity on the
canonical model $(\cal T_2)_{\mu (b)}$.  The algebraic structure is
just the natural algebraic structure on the finite automorphism group
of the stable surface.
If our covering variety happens to be disconnected, choose
a component $B''$ dominating $B_2$.  Finally,
we take $B_3 \rightarrow B$
to be the Galois normalization of $B'' \rightarrow B$, and let
$G$ be Galois group $\operatorname{Gal}(k(B_3)/k(B))$.  We let
$X_3$ denote the principal component of $X_2 \times_{B_2}B_3$.
We can represent a generic point of $X_3$ as a triple
$$(p,b,\psi)$$
where $b \in B_2$, $p \in (X_2)_b$, and $\psi:(X_2)_b \rightarrow
(\cal T_2)_{\mu(b)}$.  The birational map from $X_3$ to $Y_3=\cal T_2
\times_{\Sigma_2}B_3$ is just the evaluation map
$$(p,b,\psi) \longrightarrow (\psi(p),b,\psi)$$
Setting $B'=B_3$, we obtain the proposition that
stable surface singularities are inclusive.  $\square$
\newline \indent For the proof of the Correlation Theorem, we will
need to elaborate a bit on this situation.  Let
$\Sigma_3 \rightarrow \Sigma_1$ denote the Galois normalization of $\Sigma_1$
in the function field $k(B_3)$.  We have:
$$\begin{array}{ccc}
B_3  & \rightarrow & \Sigma_3 \\
\downarrow   &  \;   &   \downarrow\\
B & \dashrightarrow & \Sigma_1 \\
\end{array}$$
The bottom arrow is only birational.
Note that the map $B_3 \rightarrow \Sigma_2$
factors naturally through $\Sigma_3$;  this is just the Stein factorization.
We claim that $G$ acts naturally on this diagram.
$G$ consists of automorphisms of
$k(B_3)$ fixing the subfield $k(B_1)$.  Since $k(\Sigma_1)\subset k(B_1)$
these automorphisms fix $k(\Sigma_1)$ as well, and they restrict
to automorphisms of the elements of $k(B_3)$ algebraic over $k(\Sigma_1)$, i.e.
$k(\Sigma_3)$.
\newline \indent Now let ${\cal T}_3$ denote
${\cal T}_2 \times_{\Sigma_2} \Sigma_3$.  We shall show that $G$
also acts birationally and equivariantly on ${\cal T}_3 \rightarrow \Sigma_3$.
Let $s \in \Sigma_3$ be a general point, and $({\cal T}_3)_s$ the corresponding
fiber.  For $g \in G$, we need to describe the map
$$({\cal T}_3)_s \rightarrow ({\cal T}_3)_{g(s)}$$
By construction $({\cal T}_3)_s$ can be birationally identified
with some corresponding fiber
of $X_3 \rightarrow B_3$.  The action of $g$ maps this surface to
another fiber of
$X_3 \rightarrow B_3$, which in turn can be birationally identifed with
$({\cal T}_3)_{g(s)}$.
This gives a commutative diagram of varieties with (birational) G-actions
$$\begin{array}{ccc}
X_3 & \rightarrow & {\cal T}_3 \\
\downarrow & \; & \downarrow \\
B_3 & \rightarrow & \Sigma_3 \\
\end{array}$$
{}From the arguments in the previous paragraph we see that $X_3$ is
birational to $Y_3={\cal T}_3 \times_{\Sigma_3} B_3$, and that this
birational map respects the Galois action of $G$.  Taking quotients
under this action gives a dominant rational map
$$X\approx ({\cal T}_3 \times_{\Sigma_3} B_3) /G
  \longrightarrow {\cal T}_3 / G$$
This refined construction is crucial to the proof of the correlation
theorem, so we summarize it below:
\begin{cor}
Let $f:X \rightarrow B$ be a proper morphism of integral varieties.
Assume that the general fiber of $f$ is a surface of general type.
Then there exists a generically finite Galois base extension
$$B' \rightarrow B$$ with Galois group $G$, and
a finite cover of the image of $B$ in the moduli space
$$\Sigma' \rightarrow \Sigma$$
with the following properties:
\begin{enumerate}
\item {There is a tautological family of surfaces
$${\cal T}' \rightarrow {\Sigma}'$$
over $\Sigma'$.}
\item {G acts on $\Sigma'$, and this action lifts to
a $G$ equivariant rational dominant map
$$\begin{array}{ccc}
X' & \rightarrow & {\cal T}' \\
\downarrow  & \; & \downarrow \\
B' & \rightarrow & {\Sigma}' \\
\end{array}$$}
\item { The pull back of ${\cal T}'$ to $B'$ is birational to $X'$, and
the quotient of this variety under the $G$-action is birational to $X$.}
\end{enumerate}
\end{cor}

\section{Stable Surface Singularities are Negligible}
In this section, we will restrict our attention to families of
stable surfaces $f:X \rightarrow B$ over a smooth base $B$, and
their fiber products $f^n: X^n_B=X \times_B ... \times_B X \rightarrow B$.
For such families, having canonical singularities
(rational double points) is an open condition.  That is, the locus
$S \subset B$ corresponding to singularities worse than rational double
points is Zariski closed.  Here we will assume that it is a proper
subvariety of $B$.  We will prove the following:
\begin{prop}
Let $f: X \rightarrow B$ be a family of stable surfaces over
a smooth proper base $B$.
Assume the generic fiber has only canonical singularities.
Then the fiber products of this family
$$f^n: X_B^n \longrightarrow B$$
have canonical singularities.
\end{prop}
We prove this in three steps.  First, we prove a general lemma on the
singularities of fiber products.  Then we establish the result in the case
where $B=\Delta$ a complex disc.
The third step is to reduce the general case to
this special case.  For this reduction, we will utilize results
of Stevens \cite {St} on families of varieties with canonical singularities.
\newline \indent Our first lemma gives some rough information on the
singularities of fiber products:
\begin{lm}
Let
$f: X\rightarrow B$ a family of stable surfaces, such that the general
fiber is normal.
Then the $n^{th}$ fiber product
$$f^n: X_B^n \longrightarrow B$$
is a normal $\Bbb Q$-Gorenstein variety.
\end{lm}
$X_B^n$ is irreducible, because the family $X \rightarrow B$ is
flat with general fiber irreducible.  We show that $X_B^n$ is Cohen-Macaulay.
$X_B$ itself is Cohen-Macaulay, as it is a flat family of
Cohen-Macaulay varieties over a smooth base.  In particular, the dualizing
complex of the morphism $X \rightarrow B$ has only a single term, the relative
dualizing sheaf $\omega_{X/B}$.  This sheaf is flat over $B$, so the dualizing
complex of $X_B^n \rightarrow B$ is just the tensor product of the dualizing
complexes of each of the factors.  In particular, this complex has only one
term
$$\omega_{X_B^n/B}=\pi_1^* \omega_{X/B} \otimes ... \otimes
	\pi_n^*\omega_{X/B}$$
and so $X^n_B$ is Cohen-Macaulay.  Note that this implies that
$X^n_B$ satisfies Serre's condition $S_r$ for every $r>0$.
\newline \indent $X^n_B$ is reduced, because it is smooth at the generic point
and satisfies the $S_1$ condition.  We now prove that $X_B^n$ is normal.
Because $X_B^n$ satisfies the $S_2$ condition, we just need to show that
it is smooth in codimension one.  Let
$$\pi_j : X_B^n \longrightarrow X $$
be the $j^{th}$ projection map.
The singularities of $X_B^n$ are contained in the set of points
where $f^n$ fails to be a smooth morphism.  But if $f^n$ fails
to be smooth at $p$, then $f$ fails to be smooth at $\pi_j(p)$ for some $j$.
Since $f$ is smooth on a set with codimension two complement, so is
$f^n$.  Thus the singularities of $X_B^n$ are in codimension two.
\newline \indent Now we prove the $\Bbb Q$ Gorenstein assertion.
First we check that $X$ itself is $\Bbb Q$ Gorenstein, i.e.
$\omega^{[N]}_X$ is locally free for some $N$.  Since $X$ is a family
of stable surfaces, there exists an integer $N$ such that for each $b\in B$
$\omega^{[N]}_{X/B}|X_b$ is locally free.  Since $\omega_{X/B}^{[N]}$ is free
on every fiber of $X\rightarrow B$, $\omega_{X/B}^{[N]}$ is locally free.
Since $B$ is smooth, $\omega_B$ is locally free, and
$$\omega^{[N]}_X=\omega_{X/B}^{[N]} \otimes f^*\omega_B^N$$
This formula is not hard to prove.  It is
certainly true on the open set $U$ where $X \rightarrow B$ is smooth.
The complement of $U$ has codimension two by hypothesis.  Since $X$ is normal,
and both sheaves are reflexive, the formula extends to all of $X$.
For the basic properties of reflexive sheaves used here, see \cite{H}.
\newline \indent Now we prove $X^n_B$ is $\Bbb Q$ Gorenstein.
As in the previous paragraph, we have the formula
$$\omega^{[N]}_{X_B^n}=\omega^{[N]}_{X_B^n/B} \otimes {f^n}^*\omega_B^{\otimes
N}$$
so it suffices to prove that $\omega^{[N]}_{X_B^n/B}$ is locally free.
Using the general formula:
$$\omega_{X_B^n/B}=
\pi_1^*\omega_{X/B}\otimes ...\otimes \pi_n^*\omega_{X/B}\quad (*)$$
we will prove
$$\omega_{X_B^n/B}^{[N]}=
\pi_1^*\omega_{X/B}^{[N]}\otimes ...\otimes \pi_n^*\omega_{X/B}^{[N]}
\quad (**)$$
and so $\omega_{X_B^n/B}^{[N]}$ is locally free.  The left hand
side of (**) is reflexive by construction, and the right hand
side is locally free because it is the tensor product of locally free sheaves.
So we just need to prove the equivalence of
the two sides of $(**)$ on an open set with codimension two complement.
Again, we choose the open
set where $f^n$ is smooth as a morphism.  On this set, the formula
follows immediately from (*), as the dualizing sheaves are already
locally free.
This completes the proof of the lemma.  $\square$
\newline \indent Now we prove our proposition in the case where the base $B$
is one dimensional.  In this special case it takes the
following form:
\begin{prop}
Let $f: X \rightarrow \Delta $ be a family of stable surfaces over the disc.
Assume that the general fiber has only rational double points.
Then the fiber products of this family over $\Delta$ have canonical
singularities.
\end{prop}
We apply semistable reduction to the family $X \rightarrow \Delta$.
Let
$$\tilde{\Delta} \longrightarrow \Delta$$
be the ramified base change, and
$$d: Y \longrightarrow \tilde X$$
a resolution of singularities such that
all the fibers of the induced map
$$F=d \circ \tilde f : Y \longrightarrow \tilde \Delta$$
are reduced normal crossings divisors.  We have the diagram:
$$\begin{array}{ccccc}
 Y    &    \;    &    \;     &    \;    &   \;    \\
\;    & \stackrel{d}{\searrow} &    \;         &    \;
	&      \;       \\
\;    &     \;    &  \tilde X     & \rightarrow &  X   \\
\;    &    \;   &{\scriptstyle {\tilde f}}  \downarrow
	&  \;  & \downarrow  \scriptstyle{f} \\
\;  &  \;  & \tilde \Delta & \rightarrow & \Delta
\end{array}$$
Using theorem 3, we find that $\tilde X$ still
has canonical singularities.
\newline \indent
The next step is to take the nth fiber products of all the varieties in this
diagram.
We take the fiber products over the bases $\Delta $ and $\tilde \Delta$,
and we use $f^n,{\tilde f}^n$, and $d^n$ to denote the maps on the
fiber products induced by $f, {\tilde f}$, and $d$ respectively.  We have the
following diagram:
$$\begin{array}{ccrcl}
Y_{\tilde \Delta}^n &  \;  &   \;   &   \;   &  \;    \\
\;  & \stackrel{d^n}{\searrow} &   \;  &   \;   &  \;   \\
\;  &   \;  &  \tilde X_{\tilde \Delta}^n &  \longrightarrow
	&  X_{\Delta}^n  \\
\;  &  \;  & {\scriptstyle{\tilde{f}^n}} \downarrow &  \;  & \downarrow
	\scriptstyle{f^n}\\
\;  &   \;   & \tilde \Delta & \longrightarrow   & \Delta
\end{array}$$

\indent  The general lemma implies that $X_{\Delta}^n$ and
$\tilde X_{\tilde \Delta}^n$ are both $\Bbb Q$-Gorenstein and normal.
As for $Y_{\tilde \Delta}^n$,
recall that we constructed $Y$ so that its fibers over $\tilde \Delta$ have
only
reduced normal crossings.   Using an argument of Viehweg \cite {V} \S 3.6,
we see that the singularities of $Y_{\tilde \Delta}^n$ are canonical.
(Using local
analytic coordinates, we can see that the singularities are toroidal,
and so are rational.  The equations also show that the singularities
are local complete intersections, hence Gorenstein.  But rational Gorenstein
singularities are canonical).

\indent First, note that for any $M$ there is an inclusion map:
$$d_*\omega_{Y}^{M} \hookrightarrow \omega_{\tilde X}^{[M]}$$
This is because the resolution $d:Y \rightarrow {\tilde X}$ is an isomorphism
on an open set with codimension two in ${\tilde X}$,
so pluricanonical forms on $Y$
yield sections of $\omega_{\tilde X}^{[M]}$.
Moreover, since $\tilde X$ has canonical singularities
we have that this is an isomorphism for some $M$, i.e.
$$d_*\omega_{Y}^{M}=\omega_{\tilde X}^{[M]}$$
This expresses the fact that regular pluricanonical
differentials on the
smooth locus of $\tilde X$ lift to smooth differentials
on the desingularization $Y$.  We will show:
$$\omega_{{\tilde X}_{\tilde \Delta}^n}^{[M]}
	=d^n_*\omega_{Y_{\tilde \Delta}^n}^{M}  \quad (1)$$
This combined with the fact that
${\tilde X}_{\tilde \Delta}^n$ is $\Bbb Q$-Gorenstein and
$Y_{\tilde \Delta}^n$ is canonical implies
that ${\tilde X}_{\tilde \Delta}^n$ is canonical as well.
\newline \indent We prove that $(1)$ holds.  Again, we have projection
maps, which fit into a commutative diagram
$$\begin{array}{rcl}
Y_{\tilde \Delta}^n & \stackrel{\phi_j}{\rightarrow} & Y\\
{\scriptstyle{d^n}} \downarrow	& \;  &  \downarrow  \scriptstyle{d} \\
{\tilde X}_{\tilde \Delta}^n  &  \stackrel{\pi_j}{\rightarrow}	&
	{\tilde X}\\
\end{array}$$
An important element in the proof $(1)$ is the equation:
$$d^n_*\phi_j^*\omega_{Y/{\tilde \Delta}}=\pi_j^*d_*\omega_{Y/{\tilde \Delta}}
\quad (2)$$
For simplicity, we prove this for $j=1$.  We begin factoring
$\phi_1=q \circ p$ and $d^n=r\circ p$:
$$\begin{array}{ccl}
Y\times ...\times Y  	&      \;    &        \;\\
{\scriptstyle{p}} \downarrow  & \stackrel{\phi_1}{\searrow}    &	  \;\\
Y\times \tilde X \times ....\times \tilde X  &
	\stackrel{q}{\rightarrow}     &     Y\\
{\scriptstyle{r}} \downarrow  &    \;   &  \downarrow \scriptstyle{d}\\
{\tilde X}_{\tilde \Delta}^n &  \stackrel{\pi_1}{\rightarrow}
	& \tilde X\\
\end{array}$$
We set $p=\text{Id} \times d^{n-1}$, $q$
the projection onto the first factor, and $r=d \times \text{Id}^{n-1}$.
Note that $\pi_1$ is flat and the square part of the diagram is a flat base
change of $d$, so $\pi_1^*d_*\omega_{Y/{\tilde \Delta}}=
r_*q^*\omega_{Y/{\tilde \Delta}}$.
The projection formula tells us that $p_*p^*(q^*\omega_{Y/{\tilde \Delta}})
=p_*{\cal O}_{Y_{\tilde \Delta}^n} \otimes q^*\omega_{Y/{\tilde \Delta}}$.
Since $p$ is a birational map of normal varieties, we have
$p_* {\cal O}_{Y_{\tilde \Delta}^n}=
	{\cal O}_{Y \times {\tilde X}\times ...\times {\tilde X}}$.
Putting all this together gives
\begin{eqnarray*}
d^n_*\phi_1^*\omega_{Y/{\tilde \Delta}}
&=&(r \circ p)_*(q \circ p)^* \omega_{Y/{\tilde \Delta}} \\
&=& r_*p_*p^*q^*\omega_{Y/{\tilde \Delta}} \\
&=& r_*(p_*{\cal O}_{Y_{\tilde \Delta}^n}\otimes q^*\omega_{Y/{\tilde
\Delta}})\\
&=& r_*q^*\omega_{Y/{\tilde \Delta}}\\
&=& \pi_1^*d_* \omega_{Y/{\tilde \Delta}}
\end{eqnarray*}
This proves equation $(2)$.
\newline \indent In the course of proving lemma 1, recall that we established
the equation:
$$\omega_{\tilde X_{\Delta}^n/{\tilde \Delta}}^{[M]}=
\pi_1^*(\omega_{{\tilde X}/{\tilde \Delta}}^{[M]})\otimes ...
\otimes \pi_n^*(\omega_{{\tilde X}/{\tilde \Delta}}^{[M]})\quad (*)$$
Using this along with $(2)$ gives
\begin{eqnarray*}
\omega_{{\tilde X}_{\tilde \Delta}^n/{\tilde \Delta}}^{[M]}&=&
\pi_1^*(\omega_{{\tilde X}/{\tilde \Delta}}^{[M]})\otimes ...
\otimes \pi_n^*(\omega_{{\tilde X}/{\tilde \Delta}}^{[M]})\\
&=&\pi_1^*(d_*\omega_{Y/{\tilde \Delta}}^{ M}) \otimes ...
\otimes \pi_n^*(d_*\omega_{Y/{\tilde \Delta}}^{M})\\
&=&d^n_*\phi_1^*\omega_{Y/{\tilde \Delta}}^{M} \otimes ...
\otimes d^n_*\phi_n^*\omega_{Y/{\tilde \Delta}}^{M}\\
&=&d^n_*\omega_{Y_{\tilde \Delta}^n/{\tilde \Delta}}^{M}
\end{eqnarray*}
The last step is just the formula for the dualizing sheaf of a fiber product.
This completes the proof of equation $(1)$.  We conclude that
${\tilde X}_{\tilde \Delta}^n$ has canonical singularities.
\newline \indent Before completing the proof, we need to fix some additional
notation.
Set
$$G :=\operatorname{Gal}({\tilde \Delta}/ \Delta)$$
and $j$ to be the map
$$j: {\tilde X}_{\tilde \Delta}^n \longrightarrow X_{\Delta}^n$$
induced by the base change.  Let
$$ s: Z \longrightarrow {\tilde X}_{\tilde \Delta}^n$$
be an equivariant desingularization of ${\tilde X}_{\tilde \Delta}^n$
with respect to the Galois action of $G$ (see \cite {Hi}).  Then we write
the quotient map
$$Q:Z \longrightarrow Z/G$$
Note that $Z/G$ may be singular.  Finally, the
map from $Z$ to $X_{\Delta}^n$ is $G$-equivariant,
so it factors through $Z/G$ giving
a map
$$R: Z/G \longrightarrow X_{\Delta}^n$$
This is summarized in the following diagram:
$$\begin{array}{rcl}
\;		&	\; \atop Q		&	\;	\\
Z		&	\longrightarrow		&	Z/G	\\
{\scriptstyle s} \downarrow	& \;\atop j	& \downarrow {\scriptstyle R}\\
{\tilde X}_{\tilde \Delta}^n    &	\longrightarrow  &  X_{\Delta}^n	\\
 \downarrow	&	\;			& \downarrow	\\
{\tilde \Delta} &	\rightarrow		& \Delta	\\
\end{array}$$
\newline \indent Now we show that $X_{\Delta}^n$ has
canonical singularities.  Let
$$\alpha \in \Gamma(\omega_{X_{\Delta}^n}^{[m]})$$
be an $m$-pluricanonical form on $X_{\Delta}^n$.  We want to show that $\alpha$
is a smooth form on $X_{\Delta}^n$, i.e. for any desingularization
$V$ of $X_{\Delta}^n$, the
pull back of $\alpha$ to $V$ is regular.
It suffices to show that $R^*\alpha$
is a smooth form on $Z/G$, because a desingularization of $Z/G$ can also
serve as a desingularization of $X_{\Delta}^n$.  By an easy local computation,
$j^*\alpha$ vanishes to order $m(|G|-1)$ along the central fiber of
${\tilde X}_{\tilde \Delta}^n$.  Pulling back to the desingularization $Z$
(and increasing $m$ if necessary),
we see that $s^*j^*\alpha$ is a smooth form vanishing to order
$m(|G|-1)$ along the central fiber, because ${\tilde X}_{\tilde \Delta}^n$
has canonical
singularities.  The central fiber of $Z$ is precisely the fixed locus
under the action of $G$.  Therefore, we can apply the following lemma to
$\theta=s^*j^*\alpha$ to show that it descends to a smooth form on $Z/G$:
\begin{lm}
Let $G$ be a finite group acting on the variety $Z$.  Let $W$
denote a codimension $d$ subvariety of $Z$ fixed pointwise by a subgroup
$G_W \subset G$, and $\theta$ an invariant $m$-pluricanonical form.
If $\theta$ vanishes to order at least $m(|G_W|-d)$ at every such $W$, then
$\theta$ descends to a smooth form on $Z/G$,
\end{lm}
For a proof of this, see \cite{CHM} \S 4.2.
We use $\beta$ to denote this smooth form on $Z/G$.  Using the commutative
diagram above, one can see that $\beta = R^*\alpha$, i.e. the pluricanonical
form $\alpha$ pulls back to a smooth form on $Z/G$.
This proves that $X_{\Delta}^n$ has canonical singularities.  $\square$
\newline \indent To summarize, this proves the proposition in the
special case of a one dimensional base.  The proposition is a local
statement on the base (in the analytic topology), so if it is true for families
over a disc then it is true for general one dimension families.  We use this as
the base case for induction on the dimension of the base.  We now prove the
inductive step.
\newline \indent We will use the following result of Stevens(\cite {St}).
\begin{thm}
Let $g:V \rightarrow \Delta$ be a family of proper varieties.  Assume:
\begin{enumerate}
\item $V$ is a $\Bbb Q$-Gorenstein integral variety, and the fibers of $g$ are
integral varieties.
\item The general fibers $g^{-1}(s)$ have only canonical singularities.
\item The special fiber $g^{-1}(0)$ has log terminal singularities.
\end{enumerate}
Then $V$ has canonical singularities.
\end{thm}
We apply this theorem inductively to $V=X_B^n$ to reduce the
dimension of the base.  Note that $X_B^n$ has $\Bbb Q$-Gorenstein
singularities by the lemma proven above.
Choose a local analytic coordinate $y$ on $B$, and consider the level
surfaces
$$H_s=\{ b \in B : y(b)=s \}$$
for $s \in \Delta$.  Assume that none of the $H_s$ are contained in the locus
$S$ of surfaces with singularities worse than rational double points.
Let
$$h: X \rightarrow \Delta$$
be the map associating $X|_{H_s}=f^{-1}(H_s)$ to $s\in\ \Delta$.
We also have the corresponding map
$$g=h^n: X_B^n \rightarrow \Delta$$
For all $s$, $H_s \cap S$ is again a Zariski closed proper subset of $H_s$, so
the family
$$f^n_s: X|_{H_s}^n \rightarrow H_s$$
satisfies the hypotheses of Proposition 2.
Applying the inductive hypothesis, we find that
$$g^{-1}(s)=X_B^n|_{H_s}=X|_{H_s}^n$$
is canonical for all $s$.  Since all the fibers of $g$ are canonical, and
the total space $X_B^n$ is $\Bbb Q$-Gorenstein, we can apply Stevens'
theorem to conclude that $X_B^n$ is also canonical.  This concludes
the proof that stable surface singularities are negligible. $\square$

\section{Proof of the Correlation Theorem}
In this section, we prove the correlation theorem for surfaces of
general type (Theorem 1).
We first prove a special case where the family has maximal variation
of moduli and the singularities are not too bad.  By definition, a
family has maximal variation of moduli if there are no isotrivial
connected subfamilies through the generic point.

\begin{thm}[Correlation for Families with Maximal Variation]
Let $X\rightarrow B$ be a family of stable surfaces, with projective
integral base and smooth general fiber.  Assume that the associated map
$\phi:B \rightarrow \bar{\cal M}$ is generically finite.
Then there exists a positive integer $n$ such that
$X^n_B$ is of general type.
\end{thm}
Being of general type is a birational property, so there is no
loss of generality if we take the base $B$ to be smooth.
To show that $X_B^n$ is of general type for some large $n$, we must
verify two statements:
\begin{enumerate}
\item $X_B^n$ has canonical singularities
\item $\omega_{X_B^n}$ is big
\end{enumerate}
Note that the first statement is equivalent to saying that stable
surface singularities are negligible, which was proved in the last section.
The second statement allows us to get lots of pluricanonical differentials
on $X_B^n$, and which pull back to a desingularization of $X_B^n$.

The key to the second statement is the following theorem:
\begin{thm}
Let $f:X \rightarrow B$ be a family of surfaces, such that the general
fiber is a surface of general type.  Assume this family has maximal
variation.  Then for $m$ sufficiently large, we have that
$f_*\omega^m_{X/B}$ is big.
\end{thm}
This result is proven by Viehweg in \cite{V2} (and more generally for
arbitrary dimensional fibers by Koll\'ar in \cite{K2}).
We need the following consequence of this result:
\begin{prop}
Under the hypotheses of the theorem above, $\omega_{X^n_B}$ is big.
\end{prop}
We will show that Theorem 7 implies Proposition 4.  Here $S^{[n]}$
will denote the reflexive hull of the $n$th symmetric power of a sheaf.
To say that $f_*\omega_{X/B}^m$ is big means that for any ample line bundle
$H$ on $B$ there exists an integer $n$ such that
$$S^{[n]}(f_*\omega_{X/B}^m) \otimes H^{-1}$$
is generically globally generated, i.e. the global sections of this
sheaf generate over an open set of $B$.  It is equivalent to say that
this sheaf is generically globally generated for sufficiently large $n$.
Now let $T^{[n]}$ denote the reflexive hull of the $n$th tensor power
of a sheaf.  We
claim that for sufficiently large $n$
$$T^{[n]}(f_*\omega_{X/B}^m) \otimes H^{-1}$$
is generically globally generated.  To prove this, we need a result from
representation theory:
\begin{prop}
Let $V$ be an $r$ dimensional vector space over a field
of characteristic zero, and let $T^n(V)$ and $S^q(V)$ be the
$n$th tensor power and $q$th symmetric power representations
of $Gl(V)$ respectively, and write $t=r!$.
Then each irreducible component of $T^n(V)$ is a quotient
of a representation
$$S^{q_1}(V) \otimes ... \otimes S^{q_t}(V)$$
where $q_i \geq {n\over{t+1}}$.
\end{prop}
This result is proved in \cite{H2} for arbitrary `positive' representations
$T$ of $V$.  This gives us a map
$$\bigoplus S^{[q_1]}(f_*\omega_{X/B}^m)\otimes ... \otimes S^{[q_t]}
(f_*\omega_{X/B}^m) \otimes H^{-1} \rightarrow T^{[n]}
(f_*\omega_{X/B}^m) \otimes H^{-1}$$
which is surjective over an open set of $B$.  Let $H$ be ample and
globally generated, and choose $n$
large enough to guarantee that each of the
$S^{[q_i]}(f_*\omega_{X/B}^m) \otimes H^{-1}$ is generically globally
generated.  This guarantees that the left hand side is generically
globally generated, but then so is $T^{[n]}(f_*\omega_{X/B}^m)
\otimes H^{-1}$.
\newline \indent Now we will prove that for some large $n$ the
dualizing sheaf $\omega_{X^n_B}$ is big, i.e.
for large $m$ we have
$$h^0(X^n_B,\omega^{[m]}_{X^n_B})\approx
m^{(b+2n)}$$
where $b=\text{dim}(B)$.  We restrict ourselves to values
of $m$ for which
\begin{enumerate}
\item{$\omega^{[m]}_{X/B}$ is locally free.}
\item{$f_*\omega^{m}_{X/B}$ is big.}
\end{enumerate}
First we compute the canonical bundle of $X^n_B$:
$$\omega_{X_B^n}=\omega_{X_B^n/B} \otimes {f^n}^*\omega_B
=\pi_1^*\omega_{X/B} \otimes ....\otimes \pi_n^*\omega_{X/B}
	\otimes {f^n}^* \omega_B$$
As in lemma 1, taking $m$th powers gives
$$\omega^{[m]}_{X_B^n}=\pi_1^*\omega^{[m]}_{X/B} \otimes ...\otimes \pi_n^*
\omega^{[m]}_{X/B} \otimes {f^n}^*\omega^m_B$$
Applying $f^n_*$ to this gives
\begin{eqnarray*}
f^n_*\omega^{[m]}_{X_B^n}
&=&f^n_*(\pi_1^*\omega^{[m]}_{X/B} \otimes ... \otimes \pi_n^*
\omega^{[m]}_{X/B}) \otimes \omega^m_B \\
&=&f^n_*\pi_1^*\omega^{[m]}_{X/B} \otimes ... \otimes
f^n_*\pi_n^*\omega^{[m]}_{X/B} \otimes \omega^m_B \\
&=&T^{n}(f_*\omega^{[m]}_{X/B}) \otimes \omega^m_B\\
\end{eqnarray*}
Note this is also a reflexive sheaf.
The inclusion map $\omega^m_{X/B} \rightarrow \omega^{[m]}_{X/B}$
induces a map of reflexive sheaves
$$T^{[n]}(f_*\omega^m_{X/B}) \rightarrow T^{n}(f_*\omega^{[m]}_{X/B})$$
which is an isomorphism at the generic point of $B$.
\newline \indent Let $H$ be an invertible sheaf on
$B$ so that $H \otimes \omega_B$
is very ample.
By Viehweg's theorem and proposition 5, we can choose $n$ so that
$T^{[n]}(f_*\omega^m_{X/B}) \otimes H^{-m}$ is generically
globally generated for $m$ sufficiently large.
The computations of the last paragraph show that
$f^n_*\omega^{[m]}_{X^n_B}\otimes
(H\otimes \omega_B)^{-m}$ is also generically
globally generated for sufficiently large $m$.
In particular, as this sheaf has rank on the order of $m^{2n}$,
there at least this many global sections.
By our assumption on $H$, we have that $(H\otimes \omega_B)^
m$ has on the order of $m^b$ sections varying
horizontally along the base $B$.  Tensoring, we get that
$f^n_*\omega^{[m]}_{X^n_B}$ has on the order of $m^{2n+b}$
global sections.  Thus we conclude that
$$h^0(\omega^{[m]}_{X^n_B})\approx m^{2n+b}$$
This completes the proof of the proposition and the special
case of the Correlation theorem. $\square$

\indent Now we extend this special case to prove correlation
for arbitrary families $f:X \rightarrow B$ of surfaces of general type.
Since stable surface singularities are inclusive, after a generically
finite base extension
$B' \rightarrow B$ every family of surfaces of general type dominates a family
of stable surfaces $\psi:{\cal T}' \rightarrow {\Sigma}'$
with maximal variation:
$$\begin{array}{cccc}
X' & \rightarrow & {\cal T}' & \; \\
\downarrow  & \; & \downarrow  &  (*)\\
B' & \rightarrow & {\Sigma}' &  \;  \\
\end{array}$$
Take $n^{th}$ fiber products, where $n$ is chosen to ensure that
${{\cal T}'}_{{\Sigma}'}^n$
is of general type.  Use ${X'_{B'}}^n$ to denote the component of the fiber
product dominating $B'$.
We obtain a diagram:
$$\begin{array}{cccc}
{X'_{B'}}^n & \rightarrow & {{\cal T}'_{{\Sigma}'}}^n &    \;  \\
\downarrow  & \; & \downarrow  &  (**)   \\
B' & \rightarrow & {\Sigma}' &    \;      \\
\end{array}$$
with ${X'_{B'}}^n$ dominating ${{\cal T}'_{{\Sigma}'}}^n$,
a variety of general type.
This shows that the correlation result holds if we allow ourselves to make
a finite base change before we take the fiber products.
\newline \indent Now we show that the fiber products $X_B^N
\rightarrow B$ dominate a variety of general type without taking a base change,
provided $N$ is large enough.
We will use the  corollary at the end of \S 3 to construct our map.
This corollary allows us to assume that the base extension $B' \rightarrow B$
is Galois with Galois group $G$, and that $G$ acts birationally on the entire
diagram $(*)$.  That is, $G$ acts birationally on each of the
varieties in $(*)$, and this action commutes with the morphisms
of the diagram.  It follows
that $G$ acts naturally and birationally on $(**)$,
and taking quotients gives us
$$\begin{array}{ccc}
X_B^N  &  \rightarrow  &  {{\cal T}'_{{\Sigma}'}}^N /G    \\
\downarrow   &   \;   &   \downarrow   \\
B  &  \rightarrow   &  {\Sigma}' / G \\
\end{array}$$
Setting $W={{\cal T}'_{{\Sigma}'}}^N /G$, we obtain a rational dominant map
$$h: X_B^n \longrightarrow W$$
{}From the construction, we see that this map is generically finite when
restricted to a general fiber of $X_B^n \rightarrow B$.
\newline \indent
To conclude the proof, we
need to show that $W={{\cal T}'_{{\Sigma}'}}^N /G$ is of general type, for some
sufficiently large $N$.
Specifically, we will show that enough of the $m$ pluricanonical
differentials on $V={{\cal T}'_{{\Sigma}'}}^N$ descend to smooth differentials
on $W=V/G$ to guarantee that $W$ is a variety of general type.
First, note that $G$ acts faithfully on ${\cal T}'$, but
does not necessarily act faithfully on the base ${\Sigma}'$.
Let $G'$ be the maximal quotient of $G$ acting
faithfully on ${\Sigma}'$, and set $g=|G|$.
Let $\Phi_1 \subsetneq {\Sigma}'$ be the locus of
points of ${\Sigma}'$ with nontrivial stabilizer under the $G'$ action,
$\Phi_2 \subsetneq {\Sigma}'$ the locus on the base corresponding
to fibers of ${\cal T}' \rightarrow {\Sigma}'$
fixed pointwise by a nontrivial subgroup of $G$.  Let
$D_0$ be an effective divisor on ${\Sigma}'$ containing $\Phi_1 \cup
\Phi_2$, and $D$
the pullback of $g D_0$ to $V$.   Note that for large $N$,
the support of $D\subset V$ contains all the componenets of the fixed locus
with
codimension less than $g$.  This is because the only components of
the fixed point locus with small codimension correspond to fixed fibers of
${\cal T}'$.
\newline \indent
We repeat the proof we used in the
maximum variation case, except that we choose positive $H$ so that
$H \otimes \omega_{{\Sigma'}}(-gD_0)$ is very ample.
Again, we can choose $N$ so that
$$\psi^N_*\omega_V^{[m]}(-mD)=
T^{[N]}(\psi^N_*\omega^m_{{\cal T}'/{\Sigma}'})
\otimes \omega^m_{{\Sigma}'} (-mgD_0)$$
has $m^{(\text{dim} V)}$ sections.
This guarantees that
$$h^0(\omega_V^{[m]}(-mD))\approx m^{(\text{dim} V)}$$
In other words, there are lots of $m$ pluricanonical differentials
on $V$ vanishing to high order along subvarieties of $V$ nontrivial
stabilizer and codimension less than $g$.
We apply lemma 2 of \S 4 to conclude that these
forms descend to smooth forms on $W$, i.e.  forms on $W$
that pull back to regular forms on a desingularization
of $W$.  Therefore $W$ is of general type
and the Correlation Theorem is
proved. $\square$

\section{Consequences of the Correlation Theorem}
We give some consequences of the
Correlation Theorem.  Many of these are stated in \S 6 of \cite{CHM}.
The motivating conjectures can be found in \cite{L}.
\newline \indent Recall the statement of the Geometric Lang Conjecture:
\begin{conj} [Geometric Lang Conjecture]
If $W$ is a variety of general type, the union of all irreducible, positive
dimensional subvarieties of $W$ not of general type is a proper,
closed subvariety  $\Xi_W \subset W$.
\end{conj}
We will call $\Xi_W$ the {\em Langian exceptional locus} of $W$.
The following theorem describes how the Langian exceptional locus varies
in families, if the geometric Lang conjecture is true.
\begin{thm}
Assume the Geometric Lang Conjecture.
\newline Let $f:X \rightarrow B$ be a flat family of
surfaces in projective space, such that the general fiber is an integral
surface of general type.  Then there is a uniform bound on the degree of the
Langian exceptional locus of fibers that are of general type i.e.
$$deg(\Xi_{X_b})\leq D$$
\end{thm}
By Noetherian induction, it suffices to prove the bound on an open
subvariety of $B$.
Using the Correlation Theorem, for sufficiently high fiber products
$f^n: X_B^n \rightarrow B$ we obtain a dominant rational
map to a variety of general type
$$\psi:X_B^n \rightarrow W$$
We use $Y$ to denote the Langian exceptional locus of $W$, and
$Z_1$ its preimage in $X_B^n$.  Let $Z_2$ be the union of
all positive dimensional fibers of the map
$$X_B^n \rightarrow W \times B$$
We set $Z=Z_1 \cup Z_2$;  $Z$ is a proper subvariety of $X_B^n$.
\newline \indent Consider the projection map
$$\pi_n:X_B^n \rightarrow  X_B^{n-1}$$
with fiber $\pi_n^{-1}(p)$ isomorphic to the stable surface
$X_{f^{n-1}(p)}$.
Since $Z$ is a proper subvariety of $X_B^n$, over an open set
$U_{n-1} \subset X_B^{n-1}$ the fibers of $\pi_n$ are not contained
in $Z$.  For all $p \in U_{n-1}$, the degree of $Z$ restricted to
the fiber $\pi_n^{-1}(p)$ is bounded.  At the same time, the Langian locus
$\Xi_p$ of any of these fibers is contained in $Z$, because for any
component $C\subset \Xi_p$ either $\psi(C)\subset Y$ or $\psi(C)$
is a point.  This concludes the proof.  $\square$
\newline \indent This yields many remarkable corollaries.
One example is the following
\begin{cor}
Assume the Geometric Lang Conjecture.
There exists a constant $D$ such that the sum of the degrees of all the
rational and elliptic curves on a smooth quintic surface in
$\Bbb P^3$ is less than $D$.
In particular, there is a uniform bound on the number of rational and elliptic
curves on a quintic surface.
\end{cor}
Recently, Abramovich \cite{AV} has found another proof of these
results.

\indent Now we shall discuss some number theoretic consequences
of the Correlation Theorem.  First, recall the Weak Lang Conjecture:
\begin{conj}[Weak Lang Conjecture]
If $W$ is a variety of general type defined over a number field $K$, then
the $K$-rational points of $W$ are not Zariski dense in $W$.
\end{conj}
Assuming this conjecture, the Correlation Theorem implies the following:
\begin{thm}
Assume the Weak Lang Conjecture.
\newline Let $X\rightarrow B$ be a flat family of surfaces
in projective space defined over a number field $K$ such that the
general fiber is an integral surface of general type.
For any $b\in B(K)$ for which $X_b$ is of general type,
let $N(b)$ be the sum of the degrees of the components of
$\overline{X_b(K)}$.  Then $N(b)$
is uniformly bounded;  in particular, the number of $K$ rational points
not contained in the Langian locus is uniformly bounded.
\end{thm}
The proof of this is similar to the proof of the previous theorem.
Again we do induction on the dimension of the base $B$.
First, we shall show that the rational points of the fibers must
lie on a proper subscheme of bounded degree.
Choose an integer $n$ so that there is a dominant rational map
$$\psi:X^n_B \rightarrow W$$
to a variety of general type $W$.  Let $Y$ denote be a proper subvariety of $W$
that contains its $K$ rational points, and let $Z$ be its preimage in $X^n_B$.
All the $K$ rational points of $X^n_B$ are contained in $Z$.  We use
$$\pi_j:X_B^j \rightarrow X_B^{j-1}$$
to denote the projection morphisms.  Finally, let $Z_j$ denote the maximal
closed set in $X_B^j$ whose preimage in $X^n_B$ is $Z$, and let $U_j$ be the
complement to $Z_j$.  Note that $\pi_j^{-1}(Z_{j-1}) \subset Z_j$ by
definition and that for $u\in U_{j-1}$ we have that $\pi_j^{-1}(u) \cap Z_j$ is
a proper subvariety of $\pi_j^{-1}(u)$.  We will use $d_j$ to denote the
sum of the degrees of all the components of $Z_j \cap \pi_j^{-1}(u)$,
regardless of their dimensions, and we set
$$N=\operatorname{max}_j(d_j)$$
If all the $K$ rational points of $B$ are concentrated
along a closed subset, we are done by induction.  Otherwise, pick a general
$K$ rational point $b\in B$.  Let $j$ be the smallest integer
for which $U_j \cap X_b^j(K)$ is empty, and let $u\in U_{j-1} \cap
X_b^{j-1}(K)$.  We have that $X_b=\pi_j^{-1}(u)$ and our set-up guarantees
that $X_b(K) \subset Z_j \cap \pi_j^{-1}(u)$.  In particular, since we have
chosen everything generically, we find that $X_b(K)$ is contained in a
subscheme of degree $N$.
\newline \indent  Now we complete the proof.  We have shown that the rational
points on each fiber are concentrated along a subscheme of degree $N$.
The components of this subscheme consist of points, rational and elliptic
curves, and curves of higher genus.  The rational and elliptic curves
are contained in the Langian locus, so we
ignore them, and there are at most $N$ components of dimension zero.
Therefore, we just need the following lemma:
\begin{lm}
Assume the Weak Lang Conjecture.
Let $C$ be a (possibly singular) curve in projective space of degree $N$
defined over a number field $K$.  Assume $C$ has no rational or elliptic
components.  Then there is a uniform bound on the number of $K$ rational
points on $C$.
\end{lm}
First, because the degree is bounded there are only finitely
many possibilities for the geometric genera of the components of $C$.  By the
hypothesis, these genera are all at least two, so we can apply
the uniform boundedness results for curves in \cite{CHM}.
This completes the proof of the theorem. $\square$
\newline \indent In the corollary that follows,
quadratic points are points defined
over some degree two extension of the base field.
\begin{cor}
Assume the Weak Lang Conjecture.
Fix a number field $K$, and an integer $g>2$.  Then there is a uniform
bound on the number of quadratic points lying on
a non-hyperelliptic, non-bielliptic curve $C$ of genus $g$ defined over $K$.
\end{cor}
Note that quadratic points on $C$
correspond to $K$ rational points on its symmetric square $\operatorname
{Sym}^2(C)$.  Moreover, a hyperelliptic (respectively bielliptic) system on $C$
corresponds to a rational (respectively elliptic) curve on
$\operatorname{Sym}^2(C)$ (\cite{AH}).
In particular, the curves described in the theorem
are precisely those for which $\Xi_{\operatorname{Sym}^2(C)}=\emptyset$,
and so by the theorem $\#\operatorname{Sym}^2(C)(K)$ is finite and uniformly
bounded.

\parbox{3.8in}{Brendan Hassett\\
Department of Mathematics, Harvard University\\
1 Oxford Street,  Cambridge, MA 02138\\
{\tt hassett@@math.harvard.edu}}

\end{document}